\documentclass[11pt, reqno]{article}

\usepackage{fourier}
\usepackage{amsmath}
\usepackage{amssymb}
\usepackage{amsthm}
\usepackage{mathtools}
\usepackage[pdftex, bookmarks, final, colorlinks, linkcolor=black, citecolor=black, urlcolor=blue]{hyperref}
\usepackage{pdfsync}
\usepackage[nottoc]{tocbibind}
\usepackage{tensor}
\usepackage{geometry}

\newcommand{\closure}[1]{\overline{#1}}
\newcommand{\df}{\coloneqq}
\newcommand{\family}[3]{(#1)_{#2 \in #3}}
\newcommand{\functions}[2]{{}^{#1} #2}
\newcommand{\interior}[1]{{#1}^{\circ}}
\newcommand{\intersection}[2]{\bigcap_{#2} #1}
\newcommand{\mesh}[1]{\mathbf{M}(#1)}
\newcommand{\nhoodFilter}[1]{\mathcal{N}_{#1}}
\newcommand{\nhoodFilterSup}[2]{\mathcal{N}_{#1}^{#2}}
\newcommand{\nhoodCore}[1]{N_{#1}}
\newcommand{\nhoodCoreSup}[2]{N_{#1}^{#2}}
\newcommand{\set}[1]{\{ #1 \}}
\newcommand{\sequence}[2]{\transSequence{#1}{#2}{\omega}}
\newcommand{\spoke}[2]{S(#1, #2)}
\newcommand{\spokeSub}[3]{S_{#1}(#2, #3)}
\newcommand{\Tone}{\text{T}_1}
\newcommand{\transSequence}[3]{(#1)_{#2 < #3}}
\newcommand{\union}[2]{\bigcup_{#2} #1}

\renewcommand{\epsilon}{\varepsilon}

\DeclareMathOperator{\cf}{cf}
\DeclareMathOperator{\dom}{dom}
\DeclareMathOperator{\ran}{ran}

\newtheorem{claim}{Claim}
\newtheorem*{claim*}{Claim}
\newtheorem{corollary}{Corollary}[section]
\newtheorem{lemma}[corollary]{Lemma}
\newtheorem{question}[corollary]{Question}
\newtheorem{theorem}[corollary]{Theorem}
\theoremstyle{definition}
\newtheorem{definition}[corollary]{Definition}

\title{An internal characterisation of radiality}
\author{Robert Leek \\
University of Oxford}
\date{17th July 2014}

\begin{document}
\maketitle

\abstract{In this paper, we will investigate how radiality  occurs in topological spaces by considering neighbourhood bases generated by nests.
We will define a new subclass of radial spaces that contains LOTS, GO-spaces and spaces with well-ordered neighbourhood bases, called the independently-based spaces.
We show that first-countable spaces are precisely the independently-based, strongly Fr\'echet spaces and we give an example of a Fr\'echet-Urysohn space that is neither independently-based nor strongly Fr\'echet.
\\\\
\textbf{Keywords}: Fr\'echet-Urysohn, independently-based, nest, nest system, radial, spoke, spoke system, strongly Fr\'echet, well-based.
\\
\textbf{MSC (2010)}: 54A20, 54D99, 54F05}

\section{Introduction}

\emph{Radial} spaces were first introduced in \cite{herrlich_quotienten_1967} under the name of \emph{stark folgenbestimmt}, and were characterised in that paper as the pseudo-open images of LOTS (\emph{Linearly Ordered Topological Spaces}).
They are a natural generalisation of Fr\'echet-Urysohn spaces.
Although we have an external characterisation of these spaces via certain quotients of LOTS, the author felt it was insufficient to truly gain an appropriate understanding for these spaces.
Thus, the author wished to find an \emph{internal} characterisation that would lead to a deeper understanding of radiality.

The most common examples of radial spaces are LOTS, GO-spaces (\emph{Generalised Ordered} spaces) and spaces with well-ordered neighbourhood bases (e.g., first-countable spaces).
These can be \linebreak
viewed as having neighbourhoods generated by nests (sets linearly ordered by inclusion).
It is these spaces which we shall generalise from to create the class of \emph{independently-based} spaces.
These spaces have a clear picture of convergence and although they don't coincide with radial spaces (as we will show), a slight weakening of the definition does give such a characterisation of radiality in terms of a type of neighbourhood base generated by subspaces that have well-ordered neighbourhood bases at the point in question.
Briefly, these neighbourhood generators describe the `paths of convergence' to the point.

In the next section, we will introduce some terminology and basic results we shall need.
In section 3, we will introduce the class of spaces called the independently-based spaces that have neighbourhoods generated by nests in an appropriate way.
We will show that all GO-spaces and \emph{well-based} spaces (spaces with well-ordered neighbourhood bases) are independently-based and that independently-based spaces are radial.
Moreover, the property of being independently-based is \linebreak
hereditary.
In section 4, we will give a characterisation of radiality using a different kind of neighbourhood generator and prove that first-countable spaces are precisely the independently-based, strongly Fr\'echet spaces. To finish with, we will construct a Fr\'echet-Urysohn space that is neither independently-based nor strongly Fr\'echet.

We will not be assuming any separation axioms in this paper.
For more on radial spaces and related convergence properties, see \cite{arkhangelskii_properties_1980, nyikos_convergence_1992, tironi_pseudoradial_2006}.

\section{Preliminaries}

Below are some definitions that are required for this paper.

\begin{definition} \hfill
\begin{itemize}
\item A \emph{transfinite sequence} is a net with well-ordered domain.
We will use the notation $\transSequence{x_\alpha}{\alpha}{\lambda}$ for the transfinite sequence $f$ with domain $(\lambda, \in)$, where $\lambda$ is an ordinal and for all $\alpha < \lambda, f(\alpha) = x_\alpha$.

\item We say that a (topological) space $X$ is \emph{radial} at a point $x$ if for every subset $A$ of $X$ that contains $x$ in its closure, there is a transfinite sequence converging to $x$ whose range lies in $A$.
If a space is radial everywhere then we call it a \emph{radial} space.

\item A space $X$ is \emph{strongly Fr\'echet} at a point $x$ if for every decreasing sequence $\sequence{A_n}{n}$ of subsets of $X$ with $x \in \intersection{\closure{A_n}}{n < \omega}$, there exists a sequence $\sequence{x_n}{n}$ that converges to $x$ such that $x_n \in A_n$ for all $n < \omega$.
If a space is strongly Fr\'echet everywhere then we call it a \emph{strongly Fr\'echet} space.

\item A space $X$ is said to be \emph{well-based} at $x$ if $x$ has a neighbourhood base well-ordered by $\supseteq$.
Such a neighbourhood base is said to be \emph{well-ordered neighbourhood base} and if $X$ is well-based at every point it is called a \emph{well-based} space.

\item A transfinite sequence is said to \emph{converge strictly} to a point $x$ in a space if it converges to $x$ and $x$ is not in the closure of any of the initial segments of the transfinite sequence.

\item A \emph{nest} is a non-empty set linearly ordered by inclusion.

\item For a linearly ordered set $(X, <)$, we define its \emph{cofinality} to be the least ordinal $\alpha$ such that there exists a cofinal map $f : (\alpha, \in) \to (X, <)$.
This will be denoted by $\cf(X, <)$.

\item For a point $x$ in a space $X$, we denote its neighbourhood filter by $\nhoodFilterSup{x}{X}$, or $\nhoodFilter{x}$ when the space is unambiguous.
We define its \emph{neighbourhood core} to be the intersection of all neighbourhoods of $x$.
This will be denoted by $\nhoodCore{x}$.
Note that in a $\Tone$-space, $\nhoodCore{x} = \set{x}$.

\item A point $x$ in a space $X$ is \emph{quasi-isolated} if $\nhoodCore{x}$ is open (or equivalently, $\nhoodCore{x}$ is a neighbourhood of $x$).
\end{itemize}
\end{definition}

See \cite{engelking_general_1989} for any undefined notions.

The following lemma shows that we can impose further restrictions on the transfinite sequences that are witnesses to radiality.

\begin{lemma} \label{Injective strict sequence lemma}
Let $X$ be a space, $x \in X$ be radial and let $A \subseteq X$ be given such that $x \in \closure{A}$.
Then there exists an regular cardinal $\lambda \leq |X|$ and an injective transfinite sequence $\transSequence{x_\alpha}{\alpha}{\lambda}$ contained in $A$ that converges strictly to $x$.
\begin{proof}
If $A \cap \nhoodCore{x} \ne \emptyset$, pick $y \in A \cap \nhoodCore{x}$ and define $\lambda \df 1, x_0 \df y$. Now suppose otherwise, so by radiality there exists an ordinal $\lambda$ and a transfinite sequence $\transSequence{x_\alpha}{\alpha}{\lambda}$ that converges to $x$.
Choose $\lambda$ minimal and note that if $\lambda$ is not regular then there is an $\gamma < \lambda$ and a strictly increasing, cofinal map $g : \gamma \to \lambda$, so $\transSequence{x_{g(\alpha)}}{\alpha}{\gamma}$ converges to $x$; this is a contradiction by the minimality of $\lambda$.
Hence $\lambda$ is regular.

Let $\alpha < \lambda$ be given and suppose $x \in \closure{\set{x_\beta : \beta < \alpha}}$.
Then by the work above, there exists a regular cardinal $\kappa$ and a transfinite sequence $\transSequence{y_\beta}{\beta}{\kappa}$ contained in $\set{x_\beta : \beta < \alpha}$ that converges to $x$.
By minimality $\kappa \geq \lambda > |\alpha| \geq |\set{y_\beta : \beta < \kappa}|$, so by regularity there exists a $B \subseteq \kappa$ unbounded such that for all $\alpha, \beta \in B, y_\alpha = y_\beta$.
Thus $\transSequence{y_\alpha}{\alpha}{B}$ converges to $x$ and so $y_{\min(B)} \in A \cap \nhoodCore{x}$, which is a contradiction.
Hence $x \notin \closure{\set{x_\beta : \beta < \alpha}}$ and thus $\transSequence{x_\alpha}{\alpha}{\lambda}$ converges strictly to $x$.
Moreover, for all $\alpha < \lambda$
\begin{equation*}
\set{x_\beta : \beta < \lambda} \backslash \set{x_\beta : \beta < \alpha} \ne \emptyset.
\end{equation*}
For every $\alpha < \lambda$, define by transfinite recursion:
\begin{equation*}
g(\alpha) \df \min(\set{\beta < \lambda : \forall \gamma < \alpha, \beta > g(\gamma) \text{ and } x_\beta \ne x_{g(\gamma)}})
\end{equation*}
Then $g : \lambda \to \lambda$ is strictly increasing and $\transSequence{x_{g(\alpha)}}{\alpha}{\lambda}$ is injective and converges strictly to $x$.
Furthemore, injectivity implies that $\lambda \leq |X|$.
\end{proof}
\end{lemma}

\section{Independently-based spaces}

The inspiration for this new class of spaces comes from asking the question what makes certain spaces radial, particularly in relation to neighbourhood bases generated by families of nests.

\begin{definition}[Nest system]
For a non-empty family of nests $\mathcal{C} = \family{\mathcal{C}_i}{i}{I}$, we define its \emph{mesh} to be:
\begin{equation*}
\mesh{\mathcal{C}} \df \Bigg\{ \intersection{C_i}{i \in I} : \forall i \in I, C_i \in \mathcal{C}_i \Bigg\}
\end{equation*}
If it is understood that each of the nests consists of subsets of some fixed set $X$, then for each $i \in I$ we define the $i$-th spoke of $\mathcal{C}$ to be $\spoke{X}{\mathcal{C}}{i} \df \intersection{\bigcap \mathcal{C}_j}{j \in I \backslash \set{i}}$, where this will equal $X$ if $I = \set{i}$.
We will omit the subscript when there is no ambiguity.

Let $X$ be a space, $x \in X$ be given and let $\mathcal{C} = \family{\mathcal{C}_i}{i}{I}$ be a family of nests of neighbourhoods of $x$.
We shall call $\mathcal{C}$ a \emph{nest system} for $x$ if $\mesh{\mathcal{C}}$ is a neighbourhood base for $x$.
\end{definition}

It is known (see \cite{arkhangelskii_properties_1980}) that well-based spaces are radial since each point has a neighbourhood base which is a nest and the transfinite sequences needed are easily constructed from these.
What about spaces with neighbourhoods generated by two nests?
Consider the Tychonoff plank $X \df (\omega + 1) \times (\omega_1 + 1)$, where each factor is topologised by the the linear ordinal ordering.
Define $A \df \omega \times \omega_1$ and note that $(\omega, \omega_1) \in \closure{A}$.
However, if $\transSequence{x_\alpha}{\alpha}{\lambda}$ is a transfinite sequence contained in $A$ that converges to $(\omega, \omega_1)$ then $\transSequence{\pi_1(x_\alpha)}{\alpha}{\lambda}$ and $\transSequence{\pi_2(x_\alpha)}{\alpha}{\lambda}$ converge to $\omega$ and $\omega_1$ respectively, where $\pi_i$ is the projection map from $X$ to its $i$th-factor.
It is easily seen that such a transfinite sequence cannot exist since $\omega$ is first-countable in $\omega + 1$ and $\omega_1$ is a \emph{$p$-point} (every countable intersection of neighbourhoods is a neighbourhood) in $\omega_1 + 1$.
Thus $X$ is not radial at $(\omega, \omega_1)$, but it still has a neighbourhood base generated by two nests.

However, every LOTS is radial and each point in a LOTS has a nest system consisting of at most two nests.
The main difference is that each neighbourhood $(a, b)$ of $x$ in a LOTS can be split into two parts, $(a, x]$ and $[x, b)$, each acting independently of the other in the sense that $\closure{A} = \closure{A \cap (a, x]} \cup \closure{A \cap [x, b)}$ and both $(a, x]$ and $[x, b)$ are well-based at $x$.
LOTS are the inspiration for the definition of this new class of spaces:

\begin{definition}[Independently-based]
Let $X$ be a space, $x \in X$ and let $\mathcal{C} = \family{\mathcal{C}_i}{i}{I}$ be a nest system for $x$.
We say that $\mathcal{C}$ is an \emph{independent} nest system for $x$ if it also satisfies the following condition:
\begin{equation*}
\forall C \in \prod_{i \in I} \mathcal{C}_i, \bigcap \ran(C) = \union{(C(i) \cap \spokeSub{\mathcal{C}}{i})}{i \in I}
\end{equation*}
A space will be called \emph{independently-based} if each of its points has an independent nest system.
\end{definition}

To justify this definition and the motivation preceding it, we will now show that well-based and GO-spaces are independently-based, which in turn are radial.
We will first need the following lemma:

\begin{lemma}
Let $X$ be a space, $Y \subseteq X$ and let $y \in Y$ have an ING with respect to $X$.
Then $y$ has an independent nest system with respect to $Y$.
\begin{proof}
Pick an independent nest system $\mathcal{C} = \family{\mathcal{C}_i}{i}{I}$ for $y$ with respect to $X$.
For every $i \in I$, define $\mathcal{D}_i \df \set{C \cap Y : C \in \mathcal{C}_i}$.
Then $\mathcal{D} \df \family{\mathcal{D}_i}{i}{I}$ is a non-empty family of nests of $Y$-neighbourhoods of $y$ and $\mesh{\mathcal{D}} = \set{U \cap Y : U \in \mesh{\mathcal{C}}}$ is a neighbourhood base for $y$ with respect to $Y$, so $\mathcal{D}$ is a nest system for $y$ with respect to $Y$.
Furthermore, let $D \in \prod_{i \in I} \mathcal{D}_i$ be given, so there exists a $C \in \prod_{i \in I} \mathcal{C}_i$ such that $D(i) = C(i) \cap Y$ for every $i \in I$.
Then:
\begin{equation*}
\bigcap \ran(D) = \left( \bigcap \ran(C) \right) \cap Y = \Bigg( \union{(C(i) \cap \spoke{X}{\mathcal{C}}{i})}{i \in I} \Bigg) \cap Y = \union{(D(i) \cap \spoke{Y}{\mathcal{D}}{i})}{i \in I}
\end{equation*}
Thus $\mathcal{D}$ is independent.
\end{proof}
\end{lemma}

\begin{corollary} \label{Independently-based hereditary}
The property of being independently-based is hereditary.
\end{corollary}

\begin{theorem}
Every well-based space is independently-based.
\begin{proof}
Let $X$ be a well-based space and let $x \in X$ be given.
Then there exists a well-ordered neighbourhood base $\mathcal{B}_0$ for $x$.
Thus $\mathcal{B} = \family{\mathcal{B}_i}{i}{1}$ is obviously an independent nest system for $x$, so $X$ is independently-based.
\end{proof}
\end{theorem}

\begin{theorem}
Every LOTS is independently-based.
\begin{proof}
Let $(X, <)$ be a linearly ordered set that is unbounded above and below and endow $X$ with the order topology inherited from $<$.
Let $x \in X$ be given and define:
\begin{align*}
\mathcal{C}_0 \df \set{(y, \infty) : y < x}  && \mathcal{C}_1 \df \set{(-\infty, y) : y > x}
\end{align*}
Then $\mathcal{C} = \family{\mathcal{C}_i}{i}{2}$ is a family of nests of neighbourhoods of $x$.
Note that:
\begin{gather*}
\mesh{\mathcal{C}} = \set{(y, z) : y < x < z} \\
\spokeSub{\mathcal{C}}{0} = (-\infty, x] \\
\spokeSub{\mathcal{C}}{1} = [x, \infty)
\end{gather*}
Hence $\mathcal{C}$ is a nest system for $x$.
Moreover, for every $y < x$ and $z > x$:
\begin{equation*}
(y, \infty) \cap (-\infty, z) = (y, z) = (y, x] \cup [x, z) = ((y, \infty) \cap \spoke{\mathcal{C}}{0}) \cup ((-\infty, z) \cap \spokeSub{\mathcal{C}}{1})
\end{equation*}
Therefore $\mathcal{C}$ is independent and hence $X$ is independently-based.

Now suppose that $(X, <)$ is not unbounded above and below and extend it to a linearly ordered set $(Y, <)$ that is unbounded above and below; this can be achieved by considering the following linear order sum, where $\omega^*$ is the reverse order of $\omega$ with the ordinal order:
\begin{equation*}
Y \df \omega^* + (X, <) + \omega
\end{equation*}
Then by the work above, $Y$ is independently-based when endowed with the order topology.
Thus the subspace topology on $X$ is also independently-based by Corollary \ref{Independently-based hereditary}.
However, this subspace topology coincides with the original order topology, so it follows that $X$ is independently-based.
\end{proof}
\end{theorem}

\begin{corollary}
Every GO-space is independently-based.
\begin{proof}
This follows from the previous theorem and Corollary \ref{Independently-based hereditary}.
\end{proof}
\end{corollary}

Now we will show that every independently-based space is radial.
The following lemma on spokes is needed first.

\begin{lemma} \label{Spoke intersection}
Let $X$ be a space, $x \in X$ and let $\mathcal{C} = \family{\mathcal{C}_i}{i}{I}$ be a nest system for $X$.
Then for all distinct $i, j \in I, \spokeSub{\mathcal{C}}{i} \cap \spokeSub{\mathcal{C}}{j} = \nhoodCore{x}$.
\begin{proof}
Let $i, j \in I$ be distinct.
Then since $\mesh{\mathcal{C}}$ is a neighbourhood base for $x$:
\begin{align*}
\spoke{\mathcal{C}}{i} \cap \spoke{\mathcal{C}}{j} &= \Bigg( \intersection{\left( \bigcap \mathcal{C}_k \right)}{k \in I \backslash \set{i}} \Bigg) \cap \Bigg( \intersection{\left( \bigcap \mathcal{C}_l \right)}{l \in I \backslash \set{j}} \Bigg) \\
&= \intersection{\left( \bigcap \mathcal{C}_k \right)}{k \in I} \\
&= \nhoodCore{x} \qedhere
\end{align*}
\end{proof}
\end{lemma}

\begin{theorem}
Every independently-based space is radial.
\begin{proof}
Let $X$ be an independently-based space and let $A \subseteq X, x \in \closure{A}$ be given.
Then there exists a $\mathcal{C} = \family{\mathcal{C}_i}{i}{I}$ an independent nest system for $x$.
Assume for every $i \in I$, there exists a $C_i \in \mathcal{C}_i$ such that $C_i \cap \spokeSub{\mathcal{C}}{i} \cap A = \emptyset$.
Then:
\begin{equation*}
\Bigg( \intersection{C_i}{i \in I} \Bigg) \cap A = \Bigg( \union{(C_i \cap \spokeSub{\mathcal{C}}{i})}{i \in I} \Bigg) \cap A = \emptyset
\end{equation*}
This is a contradiction as $\intersection{C_i}{i \in I}$ is a neighbourhood of $x$.
Therefore there exists an $i \in I$ such that for all $C \in \mathcal{C}_i, C \cap \spokeSub{\mathcal{C}}{i} \cap A \ne \emptyset$.
Define $\lambda \df \cf(\mathcal{C}_i, \supseteq)$, so there exists a $C : \lambda \to \mathcal{C}_i$ strictly increasing with cofinal range in $(\mathcal{C}_i, \supseteq)$.
Then for all $\alpha < \lambda$, there exists an $x_\alpha \in C(\alpha) \cap \spokeSub{\mathcal{C}}{i} \cap A$.
Let $C' \in \prod_{j \in I} \mathcal{C}_j$ be given.
Then there exists an $\alpha < \lambda$ such that $C(\alpha) \subseteq C'(i)$.
Hence for all $\beta \in \lambda \backslash \alpha$:
\begin{equation*}
x_\beta \in C(\beta) \cap \spokeSub{\mathcal{C}}{i} \cap A \subseteq C(\alpha) \cap \spokeSub{\mathcal{C}}{i} \cap A \subseteq \bigcap \ran(C')
\end{equation*}
Therefore $\transSequence{x_\alpha}{\alpha}{\lambda} \to x$ and hence $X$ is radial.
\end{proof}
\end{theorem}

An alternative characterisation of independently-based spaces is given by glueing together neighbourhoods from subspaces.

\begin{definition}
If $X$ is a topological space with a point $x$ and subspace $S$ such that $\nhoodCore{x} \subseteq S$ and $S$ is well-based at $x$.
Then we say that $S$ is a \emph{spoke} of $x$ in $X$.
[Note that the spokes of a nest system for a point are indeed spokes for that point.]

Let $X$ be a space, $x \in X$ and let $\mathcal{S} = \family{S_i}{i}{I}$ be a non-empty family of spokes of $x$.
Then we say that $\mathcal{S}$ is a \emph{spoke system} for $x$ if $\set{\union{U_i}{i \in I} : \forall i \in I, U_i \in \nhoodFilterSup{x}{S_i}}$ is a neighbourhood base for $x$.
Furthermore, if for all distinct $i, j \in I, S_i \cap S_j = \nhoodCore{x}$, then we say that $\mathcal{S}$ is \emph{independent}.
\end{definition}

We can easily translate back and forth between independent systems of nests and spokes:

\begin{theorem} \label{Equiv ING forms}
Let $X$ be a space and let $x \in X$ be given.
\begin{enumerate}
\item Let $\mathcal{C} = \family{\mathcal{C}_i}{i}{I}$ be an independent nest system for $x$.
Then $\mathcal{S} \df \family{\spokeSub{\mathcal{C}}{i}}{i}{I}$ is an independent spoke system for $x$.

\item Let $\mathcal{S} = \family{S_i}{i}{I}$ be an independent spoke system for $x$, so for all $i \in I$ there exists a well-ordered neighbourhood base $\mathcal{B}_i$ for $x$ with respect to $S_i$.
For all $i \in I$, define:
\begin{equation*}
\mathcal{C}_i \df \Bigg\{ B \cup \union{S_j}{j \in I \backslash \set{i}} : B \in \mathcal{B}_i \Bigg\}
\end{equation*}
Then $\mathcal{C} \df \family{\mathcal{C}_i}{i}{I}$ is an independent nest system for $x$.
\end{enumerate}
\begin{proof} \hfill
\begin{enumerate}
\item For each $i \in I$, define $\lambda_i \df \cf(\mathcal{C}_i, \supseteq)$, so there exists an $F_i : (\lambda_i, <) \to (\mathcal{C}_i, \supseteq)$ strictly increasing and cofinal.
Define $\mathcal{B}_i \df \set{F_i(\alpha) \cap \spokeSub{\mathcal{C}}{i} : \alpha < \lambda_i}$.
Then $\mathcal{B}_i$ is a well-ordered neighbourhood base for $x$ with respect to $\spokeSub{\mathcal{C}}{i}$. Let $C \in \prod_{i \in I} \mathcal{C}_i$ be given and for every $i \in I$, choose $\alpha_i < \lambda_i$ such that $F_i(\alpha_i) \subseteq C(i)$.
Then:
\begin{equation*}
\bigcap \ran(C) \supseteq \intersection{F_i(\alpha_i)}{i \in I} = \union{(F_i(\alpha_i) \cap \spokeSub{\mathcal{C}}{i})}{i \in I}
\end{equation*}
Therefore $\set{\union{B_i}{i \in I} : \forall i \in I, B_i \in \mathcal{B}_i}$ is a neighbourhood base for $x$ in $X$.
Hence by Lemma \ref{Spoke intersection}, $\mathcal{S}$ is an independent spoke system for $x$.

\item If $I$ is a singleton then trivially $\mathcal{C} = \family{\mathcal{C}_i}{i}{I}$ is an ING for $x$.
Suppose not and for each $i \in I$, choose $B_i \in \mathcal{B}_i$.
Then for all $k \in I$:
\begin{equation*}
S_k \cap \intersection{\Bigg( B_i \cup \union{S_j}{j \in I \backslash \set{i}} \Bigg)}{i \in I} = B_k \cap S_k = B_k \Rightarrow \intersection{\Bigg( B_i \cup \union{S_j}{j \in I \backslash \set{i}} \Bigg)}{i \in I} = \union{B_i}{i \in I}
\end{equation*}
Hence $\mathcal{C}$ is a nest system for $x$.
Note that for all $i \in I$:
\begin{align*}
\spokeSub{\mathcal{C}}{i} &= \intersection{\Bigg( \bigcap \mathcal{C}_j \Bigg)}{j \in I \backslash \set{i}} = \intersection{\Bigg( \Bigg( \bigcap \mathcal{B}_j \Bigg)}{j \in I \backslash \set{i}} \cup \union{S_k}{k \in I \backslash \set{j}} \Bigg) = \intersection{\Bigg( \Bigg( \union{S_k}{k \in I} \Bigg)}{j \in I \backslash \set{i}} \backslash (S_j \backslash \nhoodCore{x}) \Bigg) \\
&= \Bigg( \union{S_k}{k \in I} \Bigg) \backslash \Bigg( \Bigg( \union{S_j}{j \in I \backslash \set{i}} \Bigg) \backslash \nhoodCore{x} \Bigg) = S_i
\end{align*}
Choose for every $i \in I, B_i \in \mathcal{B}_i$.
Then:
\begin{equation*}
\intersection{\Bigg( B_i \cup \union{S_j}{j \in I \backslash \set{i}} \Bigg)}{i \in I} = \union{B_i}{i \in I} = \union{\Bigg( \Bigg( B_i \cup \union{S_j}{j \in I \backslash \set{i}} \Bigg) \cap S_i \Bigg)}{i \in I} = \union{\Bigg( \Bigg( B_i \cup \union{S_j}{j \in I \backslash \set{i}} \Bigg) \cap \spokeSub{\mathcal{C}}{i} \Bigg)}{i \in I}
\end{equation*}
Therefore $\mathcal{C}$ is independent. \qedhere
\end{enumerate}
\end{proof}
\end{theorem}

\section{Neighbourhood characterisation of radiality}

The obvious question to ask is whether every radial space is independently-based.
If we weaken the condition for independence slightly, then we do indeed have a characterisation of radiality.

\begin{theorem} \label{Radial generator characterisation}
Let $X$ be a space and let $x \in X$ be given.
Then the following are equivalent:
\begin{enumerate}
\item $X$ is radial at $x$.
\item $X$ has a spoke system for $x$.
\item $X$ has a spoke system $\family{S_i}{i}{I}$ for $x$ such that for all distinct $i, j \in I, x \notin \closure{(S_i \cap S_j) \backslash \nhoodCore{x}}$.
\end{enumerate}
\begin{proof}
Assume $X$ has a spoke system $\family{S_i}{i}{I}$ for $x$, so for all $i \in I$ there is a well-ordered neighbourhood base $\mathcal{B}_i$ at $x$ with respect to $S_i$.
Let $A \subseteq X$ be given such that $x \in \closure{A}$ and suppose for all $i \in I$, there exists a $B_i \in \mathcal{B}_i$ such that $B_i \cap A = \emptyset$.
Then $\left(\union{B_i}{i \in I} \right) \cap A = \emptyset$, which is a contradiction.
So there is an $i \in I$ such that $x \in \closure{A \cap S_i}^{S_i}$.
As $x$ is well-based in $S_i$, it follows that there exists a transfinite sequence in $A \cap S_i$ converging to $x$.
Therefore $X$ is radial at $x$.

Now suppose $X$ is radial at $x$.
For clarity, the rest of this proof will use a function notation (e.g. $f, g, h$, etc.) rather than a sequence notation (e.g. $\transSequence{x_\alpha}{\alpha}{\lambda}$) for transfinite sequences.
If $x$ is quasi-isolated then define $I \df \set{0}, Y_0 \df \nhoodCore{x}$.
Suppose otherwise and define:
\begin{gather*}
\mathcal{T} \df \big\{f \in \bigcup \set{\functions{\lambda}{(X \backslash \nhoodCore{x})} : \lambda \text{ is a regular, non-zero cardinal and } \lambda \leq |X|} : f \to x \text{ strictly and } f \text{ is injective} \big\} \\
\mathcal{A} \df \set{\mathcal{F} \subseteq \mathcal{T} : \forall f, g \in \mathcal{F} \text{ distinct}, f^{-1}[\ran(f) \cap \ran(g)] \text{ is bounded in } \dom(f)}
\end{gather*}
By Tukey's Lemma there exists a maximal $\mathcal{F} \in \mathcal{A}$. For all $f \in \mathcal{F}$, define $S_f \df \nhoodCore{x} \cup \ran(f)$.
\begin{claim*}
For every $f \in \mathcal{F}, \mathcal{B}_f \df \set{f[\dom(f) \backslash \alpha] \cup \nhoodCore{x} : \alpha \in \dom(f)}$ is a well-ordered neighbourhood base for $x$ in $S_f$.
\begin{proof}
Let $f \in \mathcal{F}, U \subseteq X$ be open such that $x \in U$.
As $f \to x$, there exists an $\alpha \in \dom(f)$ such that $f[\dom(f) \backslash \alpha] \subseteq U$.
By definition, $\nhoodCore{x} \subseteq U$.
Now let $\beta \in \dom(f)$ be given.
Then since $f \to x$ strictly, $x \notin \closure{f[\beta]}$, so $x \in \interior{(X \backslash f[\beta])}$ and moreover $\interior{(X \backslash f[\beta])} \cap S_f \subseteq f[\dom(f) \backslash \beta] \cup \nhoodCore{x}$.
Therefore $\mathcal{B}_f$ is a well-ordered neighbourhood base for $x$ with respect to $S_f$.
\end{proof}
\end{claim*}
Let $f, g \in \mathcal{F}$ be distinct, so there exists an $\alpha \in \dom(f)$ such that $\ran(f) \cap \ran(g) \subseteq f[\alpha]$.
Then by the previous claim there exists a neighbourhood $U$ of $x$ such that $U \cap S_f = f[\dom(f) \backslash \alpha] \cup \nhoodCore{x}$.
In particular, as $f$  is injective, $U \cap f[\alpha] = \emptyset$, so $x \notin \closure{f[\alpha]}$.
However $(S_f \cap S_g) \backslash \nhoodCore{x} = \ran(f) \cap \ran(g)$, so $x \notin \closure{(S_f \cap S_g) \backslash \nhoodCore{x}}$.

Finally, for all $f \in \mathcal{F}$, pick $\alpha_f \in \dom(f)$ and define:
\begin{equation*}
U \df \union{(f[\dom(f) \backslash \alpha_f] \cup \nhoodCore{x})}{f \in \mathcal{F}}
\end{equation*}
Suppose $U$ is not a neighbourhood of $x$.
Then by Lemma \ref{Injective strict sequence lemma} there exists a regular cardinal $\lambda \leq |X|$ and an injective transfinite sequence $f : \lambda \to X \backslash U$ that converges strictly to $x$.
Let $g \in \mathcal{F}$ be given.
Then since $\ran(f) \subseteq X \backslash U$ and for every $\beta \in \dom(g) \backslash \alpha_g, g(\beta) \in U$, it follows that $g^{-1}[\ran(f) \cap \ran(g)] \subseteq \alpha_g$ and hence $x \notin \closure{\ran(f) \cap \ran(g)}$, as $g$ converges strictly to $x$.
Thus $f^{-1}[\ran(f) \cap \ran(g)]$ is bounded in $\lambda$ also.
By maximality of $\mathcal{F}, f \in \mathcal{F}$, which is a contradiction as $f(\alpha_f) \in U$.
Therefore $U \in \nhoodFilter{x}$.
Furthermore, for all $V \in \nhoodFilter{x}$ and $h \in \mathcal{F}$, as $h \to x$, there exists a $\beta_h \in \dom(h)$ such that $h[\dom(h) \backslash \beta_h] \cup \nhoodCore{x} \subseteq V$ and so:
\begin{equation*}
\union{(h[\dom(h) \backslash \beta_h] \cup \nhoodCore{x})}{h \in \mathcal{F}} \subseteq V
\end{equation*}
Therefore $\family{S_f}{f}{\mathcal{F}}$ is a spoke system for $x$.
\end{proof}
\end{theorem}

This characterisation is the best we can achieve, since for any topological space, first-countability is equivalent to being independently-based and strongly Fr\'echet, as we now show.

\begin{lemma}[{\cite[Proposition 5.11, Theorem 5.23]{arkhangelskij_frequency_1981}}]
Every topological space is strongly Fr\'echet at each of its first-countable points.
\end{lemma}

\begin{lemma}
Let $X$ be a topological space, $x \in X$ be given such that $x$ has a finite, independent spoke system whose spokes are first-countable at $x$.
Then $X$ is first-countable at $x$.
\begin{proof}
Let $\mathcal{S} = \family{S_i}{i}{I}$ be such a spoke system.
For each $i \in I$, choose a countable neighbourhood base $\mathcal{B}_i$ for $x$ with respect to $S_i$.
Define $\mathcal{B} \df \set{\union{B_i}{i \in I} : \forall i \in I, B_i \in \mathcal{B}_i}$.
Then $\mathcal{B}$ is a countable neighbourhood base for $x$ with respect to $X$, since $I$ is finite.
Therefore $x$ is first-countable in $X$.
\end{proof}
\end{lemma}

\begin{theorem}
Let $X$ be a topological space, $x \in X$ be given such that $x$ is independently-based and strongly Fr\'echet.
Then $x$ is first-countable.
\begin{proof}
Suppose otherwise.
Then there exists an independent spoke system $\family{S_i}{i}{I}$ for $x$ and without loss of generality, assume that $x$ is not quasi-isolated in each $S_i$.
As $x$ is strongly Fr\'echet, it must be first-countable in each $S_i$, so $I$ must be infinite by the previous lemma.
Choose $J = \set{i_n : n < \omega} \subseteq I$, where each $i_n$ is distinct from the others, and define for every $n < \omega, B_n \df (\union{S_{i_m}}{m \geq n}) \backslash \nhoodCore{x}$.
As $S_i \ne \nhoodCore{x}$ and $x$ is not quasi-isolated in $S_i$, it follows that $\sequence{B_n}{n}$ is descending sequence of subsets of $X$ and $x \in \intersection{\closure{B_n}}{n < \omega}$.
Since $x$ is strongly Fr\'echet, there exists a sequence $\sequence{x_n}{n}$ that converges to $x$ such that $x_n \in B_n$ for each $n < \omega$ and thus there exists a unique $j_n \in \omega \backslash n$ such that $x_n \in S_{i_{j_n}}$.
Hence for all $m < \omega, \set{x_n : n \in \omega} \cap S_{i_m}$ is finite, so there exists $U_m$ a $S_{i_m}$-neighbourhood of $x$ missing all the $x_n$'s.
Then $U \df (\union{U_m}{m < \omega}) \cup (\union{S_i}{i \in I \backslash J})$ is a neighbourhood of $x$ and $U \cap \set{x_n : n < \omega} = \emptyset$, which is a contradiction.
Therefore $x$ is first-countable.
\end{proof}
\end{theorem}

\begin{corollary}
Let $X$ be a topological space [and let $x \in X$ be given].
Then $X$ is first-countable [at $x$] if and only if $X$ is independently-based and strongly Fr\'echet [at $x$].
\end{corollary}

An example of a strongly Fr\'echet space that isn't first-countable is the one-point compactification of an uncountable discrete space (see \cite[Example 5.12]{arkhangelskij_frequency_1981}).
By the previous corollary, such a space cannot be independently-based.
We have now answered the question posed at the beginning of this section:

\begin{theorem}
There exists a Fr\'echet-Urysohn (and thus radial) space that isn't independently-based.
\end{theorem}

The previous corollary could be viewed as saying that independently-based, Fr\'echet-Urysohn spaces and strongly Fr\'echet spaces are `orthogonal'.
However, they are not `complementary'.
We will now construct a Fr\'echet-Urysohn space that is neither independently-based nor strongly Fr\'echet.
This construction will demonstrate several techniques to reflect properties from one spoke system to another.
The following lemma and corollary will aid our proof.

\begin{lemma} \label{ING clustering sequence}
Let $X$ be a space, $x \in X$ be given and let $\family{S_i}{i}{I}$ be a spoke system for $x$.
Let $\sequence{x_n}{n}$ be contained in $X \backslash \nhoodCore{x}$ and cluster at $x$.
Then there exists an $i \in I$ such that $\sequence{x_n}{n}$ has a subsequence lieing in $S_i$ that converges to $x$.
\begin{proof}
Define $R \df \set{x_n : n \in \omega}$ and assume that for every $i \in I$, there is no subsequence of $\sequence{x_n}{n}$ contained in $S_i$ that converges to $x$.
Let $i \in I$ be given, so there exists a well-ordered neighbourhood base $\mathcal{B}_i$ for $x$ with respect to $S_i$ and suppose $x \in \closure{R \cap S_i}$.
Then by Lemma \ref{Injective strict sequence lemma}, there exists a regular cardinal $\lambda$ and a strictly increasing, cofinal chain $\transSequence{B_{i, \alpha}}{\alpha}{\lambda}$ in $(\mathcal{B}_i, \supseteq)$.
If $\lambda = 1$ then as $\nhoodCoreSup{x}{S_i} = \nhoodCoreSup{x}{X} = \nhoodCore{x} \subseteq Y_i$, it follows that $B_{i, 0} = \nhoodCore{x}$, which is a contradiction.
If $\lambda$ is uncountable then for all $n < \omega$, there exists an $\alpha_n < \lambda$ such that $x_n \notin B_{i, \alpha_n}$.
Define $\alpha \df \sup(\set{\alpha_n : n < \omega}) < \lambda$.
Then $B_{i, \alpha} \cap R = \emptyset$, which is a contradiction.
Hence $\lambda = \omega$.

Define by recursion for all $n < \omega$:
\begin{equation*}
g(n) \df \min(\set{k < \omega : k > \sup(g[n]) \text{ and } x_k \in B_{i, n}})
\end{equation*}
Then $\sequence{x_{g(n)}}{n}$ is contained in $S_i$ and converges to $x$, which is a contradiction.
Therefore there exists an $n_i < \omega$ such that $B_{i, n_i} \cap R = \emptyset$.
Hence $\left( \union{B_{i, n_i}}{i \in I} \right) \cap R = \emptyset$.
This is a contradiction as $\sequence{x_n}{n}$ clusters at $x$.
Thus there is an $i \in I$ such that $S_i$ contains a subsequence of $\sequence{x_n}{n}$ converging to $x$.
\end{proof}
\end{lemma}

\begin{corollary} \label{Reflecting neighbourhood generators}
Let $X$ be a space, $x \in X$ be given and let $\family{S_i}{i}{I}, \family{T_j}{j}{J}$ be spoke systems for $x$.
Let $i \in I$ be given and define:
\begin{equation*}
K_i \df \set{j \in J : x \in \closure{(S_i \cap T_j) \backslash \nhoodCore{x}}}
\end{equation*}
Assume:
\begin{enumerate}
\item $\family{T_j}{j}{J}$ is independent.

\item $\chi(x, S_i) = \aleph_0$.
\end{enumerate}
Then $K_i$ is finite.
\begin{proof}
Suppose $K_i$ is infinite and let $\sequence{A_n}{n}$ be a well-ordered neighbourhood base for $x$ with respect to $S_i$.
Then by assumption, for all $n < \omega$, there exists a $j_n \in J$ and $x_n \in (A_n \cap T_{j_n}) \backslash \nhoodCore{x}$ such that for all distinct $m, n < \omega, j_m \ne j_n$.
Then $\sequence{x_n}{n} \to x$ and in particular clusters at $x$, so by the previous lemma there exists a $g : \omega \to \omega$ strictly increasing and a $j \in J$ such that $\sequence{x_{g(n)}}{n} \to x$ and $\set{x_{g(n)} : n < \omega} \subseteq T_j$.
By independence, it follows that there is an $n < \omega$ such that $j = j_n$, which is a contradiction since $f^{-1}[T_j] = \set{j}$.
Therefore $K_i$ is finite.
\end{proof}
\end{corollary}

\begin{theorem}
There exists a Fr\'echet-Urysohn space that is neither independently-based nor strongly Fr\'echet.
\begin{proof}
Let $| \cdot |$ be the Euclidean norm on $\mathbb{R}^2$ and for each $x \in \mathbb{R}^2$ and $\epsilon > 0$, let $B(x, \epsilon)$ denote the open $\epsilon$-ball around $x$ given by the norm.
Denote the origin by $\mathbf{0}$ and for every $x \in \mathbb{R}^2 \backslash \set{\mathbf{0}}$, define:
\begin{gather*}
S_x \df \set{y \in \mathbb{R}^2 : |y - x| = |x|} \\
\forall n \in \mathbb{N}^+, B_{x, n} \df B(\mathbf{0}, n^{-1}) \cap S_x \\
\mathcal{B} \df \Bigg\{ \union{B_{x, n_x}}{x \in \mathbb{R}^2 \backslash \set{\mathbf{0}}} : \forall x \in \mathbb{R}^2 \backslash \set{\mathbf{0}}, n_x \in \mathbb{N}^+ \Bigg\}
\end{gather*}
Endow $\mathbb{R}^2$ with the unique topology with neighbourhood base $\mathcal{B}$ at $\mathbf{0}$ and all other points are isolated.
Denote this space by $X$.
Let $x, y \in \mathbb{R}^2 \backslash \set{\mathbf{0}}$ be distinct.
Then since distinct circles intersect in at most two points, there exists an $n_y \in \mathbb{N}^+$ such that $B_{y, n_y} \cap S_x = \set{\mathbf{0}}$.
Thus for all $n \in \mathbb{N}^+$:
\begin{equation*}
\Bigg( B_{x, n} \cup \union{B_{y, n_y}}{y \in \mathbb{R}^2 \backslash \set{\mathbf{0}, x}} \Bigg) \cap S_x = B_{x, n}
\end{equation*}
Hence $\set{B_{x, n} : n \in \mathbb{N}^+}$ is a neighbourhood base for $\mathbf{0}$ with respect to $S_x$, so $S_x$ is well-based at $\mathbf{0}$.
Therefore $\family{S_x}{x}{\mathbb{R}^2 \backslash \set{\mathbf{0}}}$ is a spoke system for $\mathbf{0}$, so by Theorem \ref{Radial generator characterisation} $X$ is radial.
Furthermore, since $\set{B(\mathbf{0}, n^{-1}) : n \in \mathbb{N}^+}$ is a pseudobase for $\mathbf{0}$ in $X$, the only transfinite sequences that converge to the origin are either eventually constant or have length $\omega$, so $X$ is Fr\'echet-Urysohn.

Note that $\nhoodCore{\mathbf{0}} = \set{\mathbf{0}}$ and for all $x \in \mathbb{R}^2 \backslash \set{\mathbf{0}}, \chi(\mathbf{0}, S_x) = \aleph_0$.

Suppose $\mathbf{0}$ has an independent spoke system $\family{T_i}{i}{I}$ and for all $i \in I$, let $\mathcal{C}_i$ be a well-ordered neighbourhood base for $\mathbf{0}$ with respect to $T_i$.
Without loss of generality, we may assume that $\mathbf{0}$ is not isolated in $T_i$ for all $i \in I$.
We first aim to show that we can refine $\family{S_x}{x}{\mathbb{R}^2 \backslash \set{\mathbf{0}}}$ to be independent.
Define:
\begin{gather*}
\forall i \in I, K_i \df \set{x \in \mathbb{R}^2 \backslash \set{\mathbf{0}} : \mathbf{0} \in \closure{(S_x \cap T_i) \backslash \set{\mathbf{0}}}} \\
\forall x \in \mathbb{R}^2 \backslash \set{\mathbf{0}}, L_x \df \set{i \in I : \mathbf{0} \in \closure{(S_x \cap T_i) \backslash \set{\mathbf{0}}}} \\
\forall F \subseteq \mathbb{R}^2 \backslash \set{\mathbf{0}}, A(F) \df \set{x \in \mathbb{R}^2 \backslash \set{\mathbf{0}} : |S_x \cap F| \geq 3}
\end{gather*}
Then by Corollary \ref{Reflecting neighbourhood generators}, $L_x$ is finite for all $x \in \mathbb{R}^2 \backslash \set{\mathbf{0}}$.
Now let $x \in \mathbb{R}^2 \backslash \set{\mathbf{0}}$ be given and for all $n < \omega$, pick $x_n \in B_{x, n} \backslash \set{\mathbf{0}}$, so $\sequence{x_n}{n} \to \mathbf{0}$.
By Lemma \ref{ING clustering sequence}, there is an $i \in I$ such that $\sequence{x_n}{n}$ has a subsequence lieing in $T_i$.
Thus $i \in L_x \ne \emptyset$.
By a similar argument, $K_i$ is also non-empty for every $i \in I$.
Also note that through any three distinct points, there is at most one circle passing through all of them, so $A(F)$ is finite for finite $F \subseteq \mathbb{R}^2 \backslash \set{0}$.

\begin{claim} \label{Claim 1}
$K_i$ is finite for every $i \in I$.
\begin{proof}
Note that for all $i \in I, \chi(\mathbf{0}, T_i) = \aleph_0$, since $X$ is Fr\'echet-Urysohn and $T_i$ is well-based at $\mathbf{0}$.
Let $i \in I$ be given and let $\sequence{C_n}{n}$ be a cofinal increasing sequence in $(\mathcal{C}_i, \supseteq)$.
Suppose $K_i$ is infinite and choose for every $n < \omega$ an $x_n \in K_i$ and $y_n \in ((C_n \cap S_{x_n}) \backslash \set{\mathbf{0}}) \backslash \union{S_x}{x \in A(\set{y_j : j < n})}$ such that for all distinct $m, k < \omega, x_m \ne x_k$.
Then $\sequence{y_n}{n} \to \mathbf{0}$ and in particular clusters at $\mathbf{0}$, so by Lemma \ref{ING clustering sequence} there exists a $g : \omega \to \omega$ strictly increasing and $y \in \mathbb{R}^2 \backslash \set{\mathbf{0}}$ such that $\sequence{y_{g(n)}}{n} \to x$ and $\set{y_{g(n)} : n < \omega} \subseteq S_y$.
Then $y \in A(\set{y_j : j < g(3)})$, so $y_{g(3)} \notin S_y$, which is a contradiction.
Therefore $K_i$ is finite.
\end{proof}
\end{claim}

\begin{claim}
For all $x \in \mathbb{R}^2 \backslash \set{\mathbf{0}}$, there exists an $\alpha_x' \in \mathbb{N}^+$ such that for all distinct $x, y \in \mathbb{R}^2 \backslash \set{\mathbf{0}}, B_{x, \alpha_x'} \cap B_{y, \alpha_y'} = \set{\mathbf{0}}$.
\begin{proof}
Let $x \in \mathbb{R}^2 \backslash \set{\mathbf{0}}$ be given and suppose that for all $n \in \mathbb{N}^+$, there exists an $x_n \in (B_{x, n} \backslash \set{\mathbf{0}}) \backslash \union{T_i}{i \in S_x}$.
Then $\sequence{x_n}{n} \to \mathbf{0}$, so by Lemma \ref{ING clustering sequence} there exists a $j \in I$ such that $\sequence{x_n}{n}$ has a subsequence lieing in $T_j$ that converges to $\mathbf{0}$.
Then $j \in L_x$, which is a contradiction.
Thus there exists an $\alpha_x \in \mathbb{N}^+$ such that $B_{x, \alpha_x} \subseteq \union{T_i}{i \in S_x}$.

Let $x, y \in \mathbb{R}^2 \backslash \set{\mathbf{0}}$ be distinct such that $\set{\mathbf{0}} \ne B_{x, \alpha_x} \cap B_{y, \alpha_y} \subseteq \union{\union{(T_i \cap T_j)}{j \in S_y}}{i \in S_x}$.
Then there exists an $i \in L_x$ and $j \in L_y$ such that $T_i \cap T_j \ne \set{\mathbf{0}}$ and so $i = j \in T_x \cap T_y$.
Thus $y \in K_i \subseteq \union{K_k}{k \in L_x}$ and therefore $F_x \df \set{y \in \mathbb{R}^2 \backslash \set{\mathbf{0}, x} : B_{x, \alpha_x} \cap B_{y, \alpha_y} \ne \set{\mathbf{0}}}$ is finite.
Since distinct circles intersect in at most two points, there exists an $\alpha_x' \in \mathbb{N}^+ \backslash \alpha_x$ such that for all $y \in F_x, B_{x, \alpha_x'} \cap B_{y, \alpha_y} = \set{\mathbf{0}}$.
As $B_{x, \alpha_x'} \subseteq B_{x, \alpha_x}$, we are done.
\end{proof}
\end{claim}

Therefore $\family{B_{x, \alpha_x'}}{x}{\mathbb{R}^2 \backslash \set{\mathbf{0}}}$ is an independent spoke system or $\mathbf{0}$.
We will now show that this leads to a contradiction. Define:
\begin{equation*}
Z \df \set{x \in \mathbb{R}^2 : |x| > 1} = \union{\set{x \in \mathbb{R}^2 : |x| > 1 \text{ and } \alpha_x' = n}}{n \in \mathbb{N}^+}
\end{equation*}
By the Baire category theorem there exists an $n \in \mathbb{N}^+$ such that $Y \df \set{x \in \mathbb{R}^2 : |x| > 1 \text{ and } \alpha_x' = n}$ is not nowhere dense in the Euclidean topology.
Thus there exists an open, non-empty subset $U$ of $\mathbb{R}^2$ such that $U \cap Y$ is dense in $U$ with respect to the Euclidean topology.
Pick $x \in U$, so there exists an $\epsilon \in (0, 1)$ and a $\theta \in \mathbb{R}$ such that $B(x, \epsilon) \subseteq U$ and $x = |x| (\cos(\theta), \sin(\theta))$.
By continuity of $\sin$, there exists a $\delta \in (-\pi / 2, \pi / 2)$ such that for all $a \in (-\delta, \delta)$:
\begin{equation*}
|\sin(a)| < \frac{\epsilon}{4 n |x| (|x| + \epsilon)}
\end{equation*}
Define $L \df \set{\lambda \cdot x : \lambda \in \mathbb{R}}$, which is closed in $\mathbb{R}^2$.
Note that
\begin{equation*}
V \df \set{c (\cos(\varphi), \sin(\varphi)) : c \in (|x| + \epsilon / 2, |x| + \epsilon), \varphi \in (\theta - \delta, \theta + \delta)}
\end{equation*}
is open in $\mathbb{R}^2$ and $(U \cap V) \backslash L$ is a non-empty open subset of $U$ with respect to the Euclidean topology.
Thus there exists a $c \in (|x| + \epsilon / 2, |x| + \epsilon)$ and $\varphi \in (\theta - \delta, \theta + \delta)$ such that $y \df c(\cos(\varphi), \sin(\varphi)) \in (V \cap Y) \backslash L$.
Define $\rho : \mathbb{R}^2 \to \mathbb{R}^2, (x, y) \mapsto (-y, x)$, which is the anticlockwise rotation about the origin by $\pi / 2$ radians.
Also define:
\begin{equation*}
z \df \frac{2 (y \bullet \rho^{-1}(x))}{|x - y|^2} \rho(y - x)
\end{equation*}
By the choice of $y$, it follows that $y \bullet \rho^{-1}(x) = |x| |y| \cos(\pi / 2 \pm (\theta - \varphi))$ and $z \ne \mathbf{0}$.
Hence:
\begin{equation*}
|z| \leq \frac{2 |x| (|x| + \epsilon) |\sin(\theta - \varphi)|}{\epsilon / 2} < \frac{1}{n}
\end{equation*}
Now let $z = (z_1, z_2), x = (x_1, x_2), y = (y_1, y_2)$.
Then:
\begin{align*}
\frac{2 (y \bullet \rho^{-1}(x))}{|x - y|^2} &= \frac{2 (y_1 x_2 - y_2 x_1)}{(x_1 - y_1)^2 + (x_2 - y_2)^2} \\
z_1 &= \frac{2 (y_1 x_2 - y_2 x_1) (x_2 - y_2)}{(x_1 - y_1)^2 + (x_2 - y_2)^2} \\
z_2 &= \frac{2 (y_1 x_2 - y_2 x_1) (y_1 - x_1)}{(x_1 - y_1)^2 + (x_2 - y_2)^2}
\end{align*}
\begin{align*}
\Rightarrow z_1^2 + z_2^2 &= \frac{4 (y_1 x_2 - y_2 x_1)^2}{(x_1 - y_1)^2 + (x_2 - y_2)^2} \\
&= 2 \cdot \frac{2 (y_1 x_2 - y_2 x_1)}{(x_1 - y_1)^2 + (x_2 - y_2)^2} ((x_2 - y_2) x_1 + (y_1 - x_1) x_2) \\
&= 2 (z_1 x_1 + z_2 x_2) \\
\Rightarrow (z_1 - x_1)^2 + (z_2 - x_2)^2 &= z_1^2 + z_2^2 + x_1^2 + x_2^2 - 2 (z_1 x_1 + z_2 x_2) \\
&= x_1^2 + x_2^2
\end{align*}
Therefore $|z - x| = |x|$, so $z \in B_{x, n}$.
Since $\rho$ is an isometry:
\begin{equation*}
\frac{2 (x \bullet \rho^{-1}(y))}{|y - x|^2} \rho(x - y) = \frac{2(\rho^{-1}(x) \bullet \rho^{-2}(y))}{|x - y|^2} (-\rho(y - x)) = z
\end{equation*}
So by symmetry $z \in B_{y, n}$, which is a contradiction.
Hence $\mathbf{0}$ has no independent spoke system and therefore $X$ is not independently-based.

Finally, we'll prove that $X$ is not strongly Fr\'echet.
This proof will use similar techniques to that of Claim \ref{Claim 1}.
Choose by recursion, for every $m, n < \omega$:
\begin{equation*}
\sigma_m(n) \in B_{(m+1, 0), n+1} \backslash \big( \bigcup \big\{S_x : \tau \in \functions{m}{\omega}, x \in A(\set{\sigma_k(\tau(k)) : k < m}) \big\} \cup \set{\mathbf{0}} \big)
\end{equation*}
Thus $\sequence{\sigma_m}{m}$ is a sequence of sequences converging to $x$.
For all $m < \omega$, define:
\begin{equation*}
C_m \df \union{\ran(\sigma_{m'})}{m \leq m' < \omega}
\end{equation*}
Note that $x \in \intersection{\closure{C_m}}{m < \omega}$.
Assume $X$ is strongly Fr\'echet, so there exists a sequence $\sigma$ that converges to $x$ such that $\sigma(n) \in C_n$ for all $n < \omega$.
By Lemma \ref{ING clustering sequence}, we may assume, without loss of generality, that there exists a $z \in \mathbb{R}^2 \backslash \set{\mathbf{0}}$ such that $\ran(\sigma) \subseteq S_z$.
Note that $\set{S_{(m + 1, 0)} \backslash \set{\mathbf{0}} : m < \omega}$ is a disjoint family, so the sequences $\sequence{\sigma_m}{m}$ have disjoint ranges.
Hence for each $n < \omega$, there exists a unique $m_n < \omega$ such that $\sigma(n) \in \ran(\sigma_{m_n})$, so there is an $l_n < \omega$ such that $\sigma(n) = \sigma_{m_n}(l_n)$.
Choose $a, b, c < \omega$ such that $m_a, m_b, m_c$ are distinct and pick an $n < \omega$ such that $m_n > \max(\set{m_a, m_b, m_c})$.
Define:
\begin{equation*}
\tau : m_n \to \omega, k \mapsto \begin{cases}
l_k &\text{if } k \in \set{m_a, m_b, m_c} \\
0 &\text{otherwise}
\end{cases}
\end{equation*}
Then $\sigma(d) = \sigma_{m_d}(l_d) = \sigma_{m_d}(\tau(d))$ for $d \in \set{a, b, c}$, so $z \in A(\set{\sigma_k(\tau(k)) : k < m_n})$ and thus $\sigma(n) = \sigma_{m_n}(l_n) \notin S_z$, which is a contradiction.
Therefore $X$ is not strongly Fr\'echet.
\end{proof}
\end{theorem}

We now conclude with a few open questions concerning independently-based spaces and spoke systems.

\begin{question}
Is there a `natural' characterisation of radial spaces in terms of nests?
\end{question}

\begin{question}
Do independently-based spaces coincide with a subclass of radial spaces with stronger convergence properties?
\end{question}

\begin{question}
How do we characterise strongly Fr\'echet spaces in terms of spoke systems?
\end{question}

\bibliographystyle{alpha}
\bibliography{Bibliography}

\end{document}